\documentclass[aop,citesort,MSNbibl,dvips]{arximspdf}
\usepackage{mathbh}

\doi{10.1214/09-AOP517}
\volume{38}
\issue{4}
\pubyear{2010}
\firstpage{1583}
\lastpage{1608}

\makeatletter

\newtheorem{theorem}{Theorem}[section]
\newtheorem{lemma}[theorem]{Lemma}
\newtheorem{corollary}[theorem]{Corollary}

\newcommand{\ind}{\mathbh{1}}

\makeatother

\begin{document}
\begin{frontmatter}

\title{The growth of the infinite long-range percolation~cluster}
\runtitle{The growth of a long-range percolation cluster}

\begin{aug}
\author[A]{\fnms{Pieter} \snm{Trapman}\corref{}\thanksref{t1}\ead[label=e1]{ptrapman@math.su.se}}
\runauthor{P. Trapman}
\affiliation{Vrije Universiteit Amsterdam and University Medical Center Utrecht}
\address[A]{Department of Mathematics\\
%Vrije Universiteit Amsterdam\\
%De Boelelaan 1081a\\
%1081 HV Amsterdam\\
%The Netherlands\\
Stockholm University\\
106 91 Stockholm\\
Sweden\\
\printead{e1}} %adresu isvedimo komanda gale!
\end{aug}

\thankstext{t1}{Supported in part by The Netherlands
Organization for Scientific Research (NWO) through a grant awarded to
Ronald Meester.}

% HISTORY:
\received{\smonth{2} \syear{2009}}
\revised{\smonth{11} \syear{2009}}

% ABSTRACT
%
\begin{abstract}
We consider long-range percolation on $\mathbb{Z}^d$, where the
probability that two vertices at distance $r$ are connected by an edge
is given by $p(r) = 1-\exp[-\lambda(r)] \in(0,1)$ and the presence or
absence of different edges are independent. Here, $\lambda(r)$ is a
strictly positive, nonincreasing, regularly varying function. We
investigate the asymptotic growth of the size of the $k$-ball around
the origin, $|\mathcal{B}_k|$, that is, the number of vertices that
are within graph-distance $k$ of the origin, for $k \to\infty$, for
different $\lambda(r)$. We show that conditioned on the origin being in
the (unique) infinite cluster, nonempty classes of nonincreasing
regularly varying $\lambda(r)$ exist, for which, respectively:

$\bullet$ $|\mathcal{B}_k|^{1/k} \to\infty$ almost surely;

$\bullet$ there exist $1 <a_1 < a_2 < \infty$ such that
$ \lim_{k\to
\infty} \mathbb{P}(a_1<|\mathcal{B}_k|^{1/k}< a_2) = 1$;

$\bullet$ $|\mathcal{B}_k|^{1/k} \to1$ almost surely.

\noindent
This result can be applied to spatial SIR epidemics. In particular,
regimes are identified for which the basic reproduction number, $R_0$,
which is an important quantity for epidemics in unstructured
populations, has a useful counterpart in spatial epidemics.
\end{abstract}

% KEYWORDS
%

\end{frontmatter}

%s1 ###
\section{Introduction and results}
%s1.1 ###
\subsection{Nearest-neighbor and long-range percolation}

Ordinary or Bernoulli nearest-neighbor bond percolation models can be
used to construct undirected random graphs in which space is explicitly
incorporated. Consider an undirected ground graph
$G_{\mathrm{ground}}=(V,E)$, in which $V$ is the set of vertices and $E$
the set of edges between vertices. The random graph
$G=G(G_{\mathrm{ground}},p)$ is obtained by removing the edges in $E$
with probability $1-p$, independently of each other. In percolation
theory, properties of the remaining graph are studied. Much effort has
been devoted to understanding the dependence of $G$ on $p$ for
$G_{\mathrm{ground}}=\mathbb{L}^d := (\mathbb{Z}^d,E_{nn})$, where
$\mathbb{Z}^d$ is the $d$-dimensional cubic lattice and $E_{nn}$ is the
set of edges between nearest neighbors, that is, vertices at Euclidean
distance $1$; see \cite{Grim99} for an extensive account on
percolation on this graph.

Long-range percolation is an extension of this model: consider a
countable vertex set $V \subset\mathbb{R}^d$. Vertices at distance
$r$ (according to some norm) share an edge with probability $p(r) =1
-e^{-\lambda(r)}$, which depends only on $r$, and the presence or
absence of an edge is independent on the presence or absence of other
edges. We refer to $\lambda(r)$ as the \textit{connection function}.
Questions similar to the questions in ordinary nearest-neighbor
percolation can be asked for properties of the random graph $G=
G(V,\lambda(r))$ obtained by long-range percolation. Note that
ordinary percolation on $\mathbb{L}^d$ is a special case of long-range
percolation with $V= \mathbb{Z}^d$ and $p(r) = p \ind(r=1)$, where
$\ind$ is the indicator function and Euclidean distance has been used.

In this paper, we consider long-range percolation on $V= \mathbb{Z}^d$
and investigate properties of the $k$-ball $\mathcal{B}_k$, the set of
vertices within graph (or chemical) distance $k$ of the origin (a
definition of the graph distance is provided below). In particular, we
analyze the asymptotic behavior of the size of this $k$-ball,
$|\mathcal{B}_k|$, for $k \to\infty$. We show that there exist
nonempty regimes of nonincreasing, positive, regularly varying
connection functions for which, respectively:
\begin{itemize}
\item$|\mathcal{B}_k|^{1/k} \to\infty$ almost surely;
\item there exist $1 <a_1 < a_2 < \infty$ such that
$ \lim_{k\to\infty} \mathbb{P}(a_1<|\mathcal
{B}_k|^{1/k}< a_2) > 0$;
\item$|\mathcal{B}_k|^{1/k} \to1$ almost surely.
\end{itemize}

%s1.2 ###
\subsection{The model and notation}\label{notation}

In this paper, we will frequently use the following notation:
$\mathbb{N}$ is the set of natural numbers, including $0$, while
$\mathbb{N}_+ := \mathbb{N} \setminus\{0\}$ is the set of
strictly positive integers. Similarly, $\mathbb{R}_+ = (0,\infty)$
consists of the strictly positive real numbers. The
ceiling of a real number $x$ is defined by $\lceil x \rceil:= \min\{y
\in\mathbb{Z}; x \leq y\}$ and its floor by $\lfloor x \rfloor:=
\max\{y \in\mathbb{Z}; y \leq x\}$. For $x,y \in\mathbb{R}$, we
define $ \sum_{i=x}^{y} f(i) :=  \sum
_{i=\lfloor x \rfloor}^{\lceil y \rceil} f(i)$.
The cardinality of a set $S$ is denoted by $|S|$.\vspace*{1pt}

The probability\vspace*{1pt} space used in this paper for long-range percolation
graphs on a countable vertex set $V \subset\mathbb{R}^d$ with \textit
{connection function} $\lambda(x,y)\dvtx \mathbb{R}^d \times\mathbb
{R}^d \to(0,\infty)$ is denoted by $(\mathcal{G}_V, \mathcal{F},
\mathbb{P})$. Here, $\mathcal{G}_V$ is the set of all simple
undirected graphs with vertex set $V$, $\mathcal{F}$ is an appropriate
$\sigma$-algebra and $\mathbb{P}$ is the product measure defined by
$\mathbb{P}(\langle x,y\rangle  \in E) = p(x,y) := 1-e^{-\lambda(x,y)}$ for $x,y
\in V$, where $\langle x,y\rangle  \in E$ denotes the event that the vertices $x \in
V$ and $y \in V$ share an edge.
We say that long-range percolation system is \textit{homogeneous} if
the connection function only depends on the distance between its
arguments, that is, $\lambda(x,y) = \lambda(\|x-y\|)$. In this paper,
$\|x\|$ denotes the $L^{\infty}$-norm of $x$, that is, for $x =
(x_1,x_2, \ldots, x_d) \in\mathbb{R}^d$, $\|x\| := \max_{1 \leq i
\leq d} x_i$, and we only consider homogeneous long-range percolation
models. Our use of the $L^{\infty}$-norm is just for mathematical
convenience and using the $L^1$- or Euclidean norm would not cause
substantial changes in this paper.

We assume that $\lambda(r)$ is nonincreasing and regularly varying,
that is, $\lambda(r)$ may be written as $r^{-\beta}L(r)$ for some
$\beta\in[0,\infty)$ and $L(r)$ is slowly varying, that is, for
every $c>0$, $ \lim_{r \to\infty} L(c r)/L(r) = 1$.
%%%For technical reasons we assume that there is a $R>0$ such that
%$L(r)$ is monotone on $[R,\infty)$.
%%Formally, we define
%%\begin{equation}\label{betadef}
%%\beta:= - \lim_{r \to\infty} \frac{\log[\lambda(r)]}{
%%\end{equation}
%%which exists because $\lambda(r)$ is regularly varying \cite{Bing89}.
%For ease of further notation we define
%%\begin{equation}
%%g(r) := -\frac{\log[\lambda(r)]}{\log[r]} -\beta= -\frac{\log[L(r)]}{
%%\end{equation}
%%Note that $ \lim_{r \to\infty} g(r) = 0$.
%%%Note that $\lambda(r) - (\lambda(r))^2/2 \leq p(r) \leq\lambda(r)$
%and $\lambda(r) \approx p(r)$ for small $p(r)$.

The random long-range percolation graph is denoted by $G_V =
G_V(\lambda(r))$. With some abuse of notation, we define $G_K:=
G_{\mathbb{Z}^d \cap[\lfloor-K/2 \rfloor,\lfloor K/2 \rfloor)^d}$
for $K \in\mathbb{N}_+$ and $G := G_{\mathbb{Z}^d}$.

A path of length $n$ consists of an ordered set of edges
$(\langle v_{i-1},v_{i}\rangle )_{1 \leq i \leq n}$. Furthermore, if the vertices $\{
v_{i}\}_{0 \leq i \leq n} \in V$ are all different, then this path is
said to be \textit{self-avoiding}. Vertices $x$ and $y$ are in the same
\textit{cluster} if there exists a path from $x$ to $y$. The graph
distance or chemical distance, $D_V(x,y) = D_{G_V}(x,y)$, between $x$
and $y$ is the (random) minimum length of a path from $x$ to $y$ in
$G_V$. If $x$ and $y$ are not in the same cluster, then $D_V(x,y) =
\infty$. Furthermore, we set $D_V(x,x)=0$. We use $D(x,y)$ for
$D_{\mathbb{Z}^d}(x,y)$ and for $K \in\mathbb{N}_+$, we define
$D_K(x,y) := D_{\mathbb{Z}^d \cap[\lfloor-K/2 \rfloor,\lfloor K/2
\rfloor)^d}(x,y)$.\vspace*{1pt}

If the probability that the origin is contained in an infinite cluster
(a cluster containing infinitely many vertices) of $G$ is positive,
then the long-range percolation system is said to be \textit
{percolating}. If a homogeneous long-range percolation system is
percolating, Kolmogorov's zero--one law (see, e.g., \cite{Grim92},
page 290) gives that $G$ almost surely contains at least one
infinite cluster, while Theorem 0 of \cite{Gand92} (see also \cite
{Berg02}, Theorem 1.3) implies that, under mild conditions, the
infinite cluster is a.s. unique. These mild conditions are satisfied
for homogeneous long-range percolation models on $\mathbb{Z}^d$, for
which $\lambda(r)$ is nonincreasing. This unique infinite cluster is
denoted by $\mathcal{C}_{\infty}$. Throughout, we will only consider
percolating systems.

For $x \in\mathbb{Z}^d$, the set $\mathcal{B}_k(x)$ is defined by
$\mathcal{B}_k(x):= \{y \in\mathbb{Z}^d; D(x,y) \leq k\}$ and
$\mathcal{B}_k :=\mathcal{B}_k(0)$. We define (as in \cite{Trap07}):
%
%e2 ###
%e1 ###
%
\begin{eqnarray}
\label{Rsteronder}
\underline{R}{}_* & := & \liminf_{k \to\infty} (\mathbb{E}(|\mathcal
{B}_k|))^{1/k} ;\\
\label{Rsterboven}
\overline{R}_* & := & \limsup_{k \to\infty} (\mathbb{E}(|\mathcal
{B}_k|))^{1/k}.
\end{eqnarray}
If $\underline{R}{}_* = \overline{R}_*$, then $R_*:= \underline{R}{}_* =
\overline{R}_*$.

%s1.3 ###
\subsection{The main results}

\begin{theorem}\label{mainthm}
Consider a percolating homogeneous long-range percolation model, as
defined in Section \ref{notation}, with vertex set $\mathbb{Z}^d$ and
nonincreasing connection function $\lambda(r) = r^{-\beta} L(r)$,
where $L(r)$ is slowly varying and $\beta\in\mathbb{R}_+$.

\textup{(a)} If either $\beta< d$, or $\beta= d$ and $
\int_{1}^{\infty} L(r) r^{-1} \,dr = \infty$,
then $\mathbb{P}(\mathcal{B}_1 = \infty) =1$. In particular,
$|\mathcal{B}_k|^{1/k} = \infty$ a.s. for $k \in\mathbb{N}_+$. So,
$R_* = \infty$.

\textup{(b)} If $\beta= d$, there exists a $K>1$ such that $L(r)$ is
nonincreasing on $[K,\infty)$ and the following conditions are satisfied
%
%e4 ###
%e3 ###
%
\begin{eqnarray}
\label{sumable}
\int_{1}^{\infty} \frac{L(r)}{r} \,dr & < & \infty,\\
\label{annoyass}
-\int_{K}^{\infty}\frac{\log[L(r)]}{r (\log[r])^2}\,dx & < & \infty,
\end{eqnarray}
then there exist constants $1 < a_1 \leq a_2 < \infty$ such that
\[
\lim_{k \to\infty} \mathbb{P}(a_1 < |\mathcal{B}_k|^{1/k} < a_2|0
\in\mathcal{C}_{\infty}) =1.
\]
Furthermore, $1<\underline{R}{}_* \leq\overline{R}_* < \infty$.

\textup{(c)} If $\beta> d$, then, for $k \to\infty$, $|\mathcal
{B}_k|^{1/k} \to1$ a.s. Furthermore, $R_*=1$.
\end{theorem}

Part (a) of this theorem is almost trivial and is only stated for
reasons of completeness.
A function which satisfies all of the conditions in part (b) is
$\lambda(r) =r^{-d} (\log[1+r])^{-\gamma}$ for $\gamma>1$. Part (b)
is the main result and perhaps the most surprising result of the paper.
Part (c) is not surprising if one knows the results of~\cite
{Bisk04}. However, some work has to be done. We prove part (c) by
using the fact that $\mathbb{P}(D(0,x) \leq n)$ decreases faster than
$\|x\|^{-\beta'}$ if $\|x\| \to\infty$. This is a result of the
following, stronger, theorem, the proof of which also provides a
simplification of the proof of the main result in \cite{Bisk04}.
\begin{theorem}\label{nonexpothm}
Consider a percolating homogeneous long-range percolation model, as
defined in Section \ref{notation}, with vertex set $\mathbb{Z}^d$ and
nonincreasing connection function $\lambda(r) = r^{-\beta} L(r)$,
where $L(r)$ is slowly varying and $\beta>d$.
Let the constants $\alpha$, $\beta'$ and $\beta''$ be such that $d <
\beta' <\beta'' < \min(2d,\beta)$ and, for all $r \geq1$,
that $\lambda(r) \leq\alpha r^{-\beta''}$.
There exists a positive constant $c = c(\alpha, \beta'', \beta')$
such that for $\gamma:= \frac{\log(2d/\beta')}{\log2} < 1$, $K(n)
:= \exp[c n^{\gamma}] +1$, all $n \in\mathbb{N}$ and all $x \in\{x
\in\mathbb{Z}^d; \|x\| > K(n)\}$, it holds that
%
%e5 ###
%
\begin{equation}
\mathbb{P}\bigl(D(0,x) \leq n\bigr) \leq[K(n)]^{\beta'} \|x\|^{-\beta'}.
\end{equation}
\end{theorem}

%s1.4 ###
\subsection{Motivation from epidemiology}\label{motivation}

We consider an SIR (Susceptible $\to$ Infectious $\to$ Recovered)
epidemic with a fixed infectious period (which, without loss of
generality, will be taken to be of length 1) in a homogeneous, randomly
mixing population of size $n$. In this model, pairs of individuals
contact each other according to independent Poisson processes with rate
$\lambda/n$. If an infectious individual contacts a susceptible one,
the latter becomes infectious as well. An infectious individual stays
infectious for one time unit and then recovers and stays immune
forever. Usually, it is assumed that there is initially one infectious
individual, with a remaining infectious period of one time unit, and
all other individuals are initially susceptible.

The basic reproduction number, $R_0$ of an SIR epidemic process in a
large, homogeneous, randomly mixing population of size $n$ is defined
as the expected number of individuals infected by a single infectious
individual in a further susceptible population \cite{Diek00}. To
proceed, we define $X^n_0$ as the set of initially infected individuals
in a population of size $n$. These individuals are said to enter
$X^n_0$ at time $0$. For $k \in\mathbb{N}$, an individual not in
$\bigcup_{j=0}^k X^n_j$ enters $X^n_{k+1}$ at the first instance it has
contact with an individual which itself entered $X^n_k$ at most one
time unit ago. We define $\mathcal{B}^n_k := \bigcup_{j=0}^k X^n_j$.
Note that the actual chain of infections that has caused the
infectiousness of an individual in $X^n_k$ might be longer than length
$k$ because it is possible that the time needed to traverse this longer
infection chain is less than the time needed to traverse the chain of
$k$ contacts that caused the individual to be in~$X^n_k$.

It has long been known (see, e.g., \cite{Ball95}) that in randomly
mixing populations, SIR epidemics can be coupled to branching
processes, in the sense that we can simultaneously define a
Galton--Watson process $\{Z_k\}_{k \in\mathbb{N}}$ (for a definition,
see \cite{Jage75}) and an epidemic processes $\{|X^n_k|\}_{k \in
\mathbb{N}}$, for all $n \in\mathbb{N}$ on one probability space,
such that for every $k \in\mathbb{N}$ and as $n \to\infty$,
$\mathbb{P}(|X_k^n| \to Z_k) =1$. In this approximation, $R_0$
corresponds to the offspring mean $m :=  \lim_{n \to
\infty} \mathbb{E}(Z_1|Z_0=1)$ of the Galton--Watson process. From
the theory of branching processes, we know that under mild conditions,
$m>1$ implies that $m^{-k} \sum_{i=0}^k Z_i$ converges a.s. to an
a.s. finite random variable which is strictly positive with nonzero probability.
By the relationship between $R_0$ and the offspring mean $m$, we deduce
that if $R_0>1$ in large populations, then the expectation $\mathbb
{E}(|\mathcal{B}^n_k|)$ will initially grow exponentially in $k$ (with
base $R_0$) and $|\mathcal{B}^n_k|$ will also grow exponentially (with
base $R_0$) with positive probability \cite{Jage75}. In particular, it
holds that
%
%e6 ###
%
\begin{equation}\label{defR0}
\lim_{k \to\infty} \lim_{n \to\infty} (\mathbb{E}(|\mathcal
{B}^n_k|))^{1/k} = \max(R_0,1).
\end{equation}

In this paper, we investigate whether it is possible to define a
quantity with similar properties as $R_0$ for spatial epidemics.

Assume that the individuals in the population are located at $\mathbb
{Z}^d$ and that the epidemic starts with one infectious individual at
the origin and all other individuals are initially susceptible. A pair
of individuals at $L^{\infty}$-distance $r$ will make contacts
according to independent Poisson processes with rate $\lambda(r)$. The
Poisson processes governing the contacts are independent. The
probability that an infectious individual makes at least one contact
with a given individual at distance $r$ during its infectious period is
given by $p(r)= 1- e^{-\lambda(r)}$. For this spatial epidemic, let
$X_k$ be defined as $X_k^n$ is defined above. It is easy to see that
the law of $\bigcup_{j=0}^k X_j$ is the same as the law of $\mathcal
{B}_k$ in the long-range percolation model with connection function
$\lambda(r)$ (see \cite{Cox88} for an exposition on this relationship
for nearest-neighbor bond percolation).

It is possible to define $R_0$ for spatial epidemics with the usual
definition $R_0 = \mathbb{E}(|X_1|||X_0|=1)$. However, this definition
is of no practical use because there is no reason to assume that
$\mathbb{E}(|X_1|||X_0|=1) =1$ is a threshold above which a large
epidemic is possible and below which it is impossible. Indeed, if $p(x)
= p$ for $x$ at Euclidean distance 1 from the origin and $0$ otherwise,
then it is known that on~$\mathbb{Z}^2$, $p =1/2$ is a threshold \cite
{Grim99,Kest80} which corresponds to $\mathbb{E}(|X_1|||X_0|=1) = 2$.
For more results on the growth of the nearest-neighbor bond percolation
cluster, see~\cite{Anta96}.

The definitions (\ref{Rsteronder}) and (\ref{Rsterboven}), and, if it
exists, the corresponding $R_*$, might be useful and might provide
information about the spread of the spatial epidemic. These definitions
are inspired by (\ref{defR0}).
Theorem \ref{mainthm} implies that regimes of $\lambda(r)$ exist in
which $R_* = \infty$, $R_* = 1$ and $1< \underline{R}{}_* \leq
\overline{R}_* < \infty$.

Note that only if $1< \underline{R}{}_* \leq\overline{R}_* < \infty$
will the quantities $\underline{R}{}_*$ and $\overline{R}_*$ appear to
be informative because $R_* =1$ does not even contain information
concerning whether an epidemic survives with positive probability or
not. Although, for $R_* = \infty$, the number of infected individuals
will be immense within a few generations, $R_*$ does not really reveal
anything about the asymptotic behavior of the spread.

A real-life application of long-range percolation for the spread of
epidemics can be found in \cite{Davi08}, where the spread of plague
among great gerbils in Kazakhstan is modelled using techniques from
(long-range) percolation theory. The present paper may be seen as the
mathematically rigorous counterpart of the paper by Davis et al.
\cite{Davi08}.

%s2 ###
\section{Remarks and discussion}\label{discussion}

Without costs, we could replace Theorem \ref{mainthm}
by the following, more general, but less elegant, theorem.
\begin{theorem}\label{mainthm2}
Consider a percolating, homogeneous, long-range percolation model, as
defined in Section \ref{notation}, with vertex set $\mathbb{Z}^d$ and
nonincreasing connection function $\lambda(r)$.

\textup{(a)} If $\sum_{x \in\mathbb{Z}^d} 1-e^{-\lambda(\|x\|)} =
\infty$, then $\mathbb{P}(\mathcal{B}_1 = \infty) =1$. Therefore,
$|\mathcal{B}_k|^{1/k} = \infty$ a.s. for $k \in\mathbb{N}_+$ and
$R_* = \infty$.

\textup{(b)} If $\lambda(r) > r^{-d}L'(r)$ is nonincreasing, $\sum
_{x \in\mathbb{Z}^d} 1-e^{-\lambda(\|x\|)} < \infty$ and $L'(r)$ is
positive, nonincreasing, slowly varying and satisfies
\[
-\int_{K}^{\infty}\frac{\log[L'(r)]}{r (\log[r])^2}\,dx < \infty,
\]
then there exist constants $a_1>1$ and $a_2 < \infty$ such that
\[
\lim_{k \to\infty} \mathbb{P}(a_1 < |\mathcal{B}_k|^{1/k} < a_2|0
\in\mathcal{C}_{\infty}) =1.
\]
Furthermore, $1<\underline{R}{}_* \leq\overline{R}_* < \infty$.

\textup{(c)} If $\liminf_{x \to\infty} -\log[\lambda(r)]/\log[r]
>d$, then for $k \to\infty$, $|\mathcal{B}_k|^{1/k} \to1$ a.s.
Furthermore, $R_*=1$.
\end{theorem}

Contrary to Theorem \ref{mainthm2}(b), part (b) of the above theorem
includes a class of connection functions that are constant on $[n,n+1)$
for every $n \in\mathbb{Z}$ and some other piecewise constant
connection functions [for which there exists no $K$ such that $L(r) =
r^d \lambda(r)$ is nonincreasing on $[K,\infty)$].

Condition (\ref{annoyass}) is troublesome because this
assumption means that this paper does not deal with all possible
nonincreasing, regularly varying connection functions.
An example of a function which does not satisfy (\ref{annoyass}), but
satisfies the other assumptions in Theorem \ref{mainthm}(c), is
\[
\lambda(r) = r^{-d} L(r) = r^{-d} \exp\biggl(-\frac{\log[r]}{\log[1+\log
[r]]}\biggr)  \qquad\mbox{for $r>1$.}
\]
%
%%For any $\gamma>1$, it holds that $\lambda(r) < r^{-d} (\log[r])^{-
%satisfied. On the other hand, we have
%%\begin{multline}
%%-\int_{K}^{\infty} \frac{\log[L(r)]}{r (\log[r])^2} \,dr = \int_{K}^{
%%\geq\int_{K}^{\infty} \frac{1}{(1+r) (1+\log[1+r])\log[1+
%% = \log[\log[1+\log[1+r]]]|_{r=K}^{\infty} = \infty.
%%\end{multline}
However, as stated above, for $\gamma> 1$, functions of the form
\[
\lambda(r) = r^{-d} (\log[1+r])^{-\gamma}
\]
satisfy all of the conditions of Theorem \ref{mainthm}(b). So, the
class of functions treated in the second statement of the theorem is
not empty.
We do not know whether Theorem \ref{mainthm}(b) still holds without
condition (\ref{annoyass}).

Theorem \ref{mainthm}(b) gives rise to some other questions,
such as:
\begin{enumerate}
\item Does $R_*$ exist for long-range percolation models
with connection functions in the regime of Theorem \ref{mainthm}(b)?
\item Does the long-range percolation graph obtained in
Theorem \ref{mainthm}(b) have a nonamenable subgraph $G'=(V',E')$?
That is, does it hold that
%
%e7 ###
%
\begin{equation}
\inf_{W \subset V';0<|W|<\infty} \frac{|\delta W|}{|W|} >0,
\end{equation}
where $\delta W$ is the set of edges in $G'$, with one end-vertex in
$W$ and one end vertex in $V' \setminus W$?
\end{enumerate}

The assumption $p(r)<1$ [i.e., $\lambda(r)<\infty$] is only
used for ease of exposition. All results in this paper are equally
valid if we relax this assumption and replace condition (\ref
{sumable}) by $ \int_{R_1+1}^{\infty} \lambda(r)
r^{d-1} \,dr = \infty$, where $R_1:= \inf\{r \in\mathbb{R};p(r) < 1\}
$. Therefore, we may allow $p(1) =1$, in order to guarantee that it is
possible to have an infinite component for any dimension $d$ and any
$\beta$ for long-range percolation on~$\mathbb{Z}^d$. Indeed, if
$d=1$ and $\beta>2$, then an infinite component only exists if $p(1) =
1$ \cite{Newm85}.

It is tempting to add the assumption $p(1)=1$ to Theorem \ref
{mainthm}(b). With that extra condition, the proofs in this paper will
become easier. However, without this extra assumption, the results of
Theorem \ref{mainthm} can be generalized to the random connection
model \cite{Mees96}, that is, long-range percolation, where the
vertex set is generated by a homogeneous Poisson point process on
$\mathbb{R}^d$. This is important in biological applications, where
exact lattice structures will not appear and models in which the
individuals/vertices are located according to a Poisson point process
might be more realistic (see, e.g., \cite{Davi08}).

Up until now, in the literature, the majority of the effort has
gone into investigating the scaling behavior of the maximum diameter of
the clusters of a homogeneous long-range percolation graph defined on
the block $V_K = \mathbb{Z}^d \cap[\lfloor-K/2 \rfloor, \lfloor K/2
\rfloor)^d$, that is, in obtaining
\[
D_K :=\max_{x,y \in V_K; D(x,y)<\infty} D_K(x,y) ;
\]
see, for example, \cite{Benj01,Bergup,Bisk04,Biskup,Copp02}. Some of
the results have been proven under the extra assumption that $p(1) =1$.

Benjamini et al. \cite{Benj04} proved that for $\lambda(r) =
r^{-\beta} L(r)$, where $\beta<d$ and $L(r)$ is slowly varying, $\lim
_{K \to\infty} D_K = \lceil d/(d - \beta) \rceil$ a.s. (see also
\cite{Benj01}). Coppersmith, Gamarnik and Sviridenko \cite{Copp02}
showed that for $\lambda(r) = \alpha r^{-d}$ and $K \to\infty$, the
quantity $D_K \log[\log[K]]/\log[K]$ is a.s. bounded away from $0$
and $\infty$.

We define $\mathcal{C}_K$ as the (random) largest cluster of the
long-range percolation graph $G_K$ (recall that $G_K := G_{\mathbb
{Z}^d \cap[\lfloor-K/2 \rfloor, \lfloor K/2 \rfloor)^d}$). In case
of a tie, $\mathcal{C}_K$ is chosen uniformly at random from the
largest clusters.
Note that if $D_K = k$ and there exists a $\rho$ such that $|\mathcal
{C}_K| > \rho K^d$ with probability tending to 1 if $K \to\infty$,
then $|\mathcal{B}_k| > \rho K^d$ with positive probability for $k \to
\infty$. So, there is an obvious relation between the diameter of a
long-range percolation cluster on $V_K$ and the rate at which $\mathcal
{B}_k$ grows. However, this relation and the results stated above do
not help us directly in obtaining Theorem \ref{mainthm}(b) and (c)
because the regime of part (b) is not even considered in the papers
cited above and the proof of the statement that $D_K \log[\log
[K]]/\log[K]$ is bounded away from 0 and $\infty$ for $\lambda(r) =
\alpha r^{-d}$ in \cite{Copp02} critically depends on the fact that
$\sum_{i\in\mathbb{N}_+} \lambda(i) i^{d-1} = \infty$. Although
results on the diameter of a long-range percolation cluster on $V_K$
may provide a lower bound for the number of vertices that are within
graph distance $k$ of the origin, they do not provide an upper bound.
So, these results are of no direct help in proving the final statement
of the theorem.

Biskup proved the following theorem (given here in our notation).
\begin{theorem}[\cite{Bisk04}]\label{biskupthm}
Consider a percolating homogeneous long-range percolation model, as
defined in Section \ref{notation}, with vertex set $\mathbb{Z}^d$ and
nonincreasing connection function $\lambda(r) =r^{-\beta}L(r)$
with\vspace*{1pt}
$\beta\in(d,2d)$ and such that $L(r)$ is positive and slowly varying.
Then, for $\Delta= \frac{\log[2]}{\log[2d/\beta]}$ and every
$\varepsilon>0$, we have
%
%e8 ###
%
\begin{equation}
\lim_{\|x\| \to\infty} \mathbb{P} \biggl(\Delta-\varepsilon\leq\frac
{\log[D(0,x)]}{\log[\log[\|x\|]]} \leq\Delta+ \varepsilon| 0,x \in
\mathcal{C}_{\infty}
 \biggr) = 1.
\end{equation}
\end{theorem}

Note that $\Delta>1$. This theorem implies that for every $\varepsilon
>0$ and every sequence of vertices $\{x_k;x_k \in\mathbb{Z}^d\}$,
\[
\lim_{k \to\infty} \ind\bigl(\|x_k\| > \exp\bigl[k^{(\Delta- \varepsilon
)^{-1}}\bigr]\bigr) \mathbb{P}\bigl(D(0,x_k) \leq k\bigr) = 0,
\]
but it does not give results on the rate at which this probability
decreases to 0. This rate is needed to prove whether or not $|\mathcal
{B}_k|^{1/k} \to1$.
Theorem \ref{nonexpothm} entails
%
%e9 ###
%
\begin{equation}
\lim_{\|x\| \to\infty} \mathbb{P} \biggl(\frac{\log[D(0,x])}{\log
[\log[\|x\|]]} \geq\Delta- \varepsilon \biggr) = 1
\end{equation}
from \cite{Bisk04} as Corollary \ref{newcorbisk}. The proofs of
Theorem \ref{nonexpothm} and Corollary \ref{newcorbisk} are shorter
and, arguably, more straightforward than the proof of the lower bound
in Theorem \ref{biskupthm} as provided in \cite{Bisk04} (cf. \cite{Biskup}).

For $\beta>2d$, Berger \cite{Bergup} proved that
%
%e10 ###
%
\begin{equation}
\liminf_{\|x\| \to\infty}  \biggl(\frac{D(0,x)}{\|x\|} \biggr) > 0
\end{equation}
almost surely. This implies that with probability 1, the growth of
$|\mathcal{B}_k|$ is of order at most $k^d$.

In a recent manuscript, Biskup \cite{Biskup} proved that if
$p(1)=1$ and $\beta$ and $\Delta$ are as in Theorem \ref{biskupthm},
then, for every $\varepsilon>0$,
\[
\lim_{L \to\infty} \mathbb{P} \bigl((\log[L])^{\Delta-\varepsilon}
\leq D_L \leq(\log[L])^{\Delta+\varepsilon}
 \bigr) = 1.
\]
Furthermore, he proved that for $p(1)=1$ and $\Lambda(r) := \mathbb
{Z}^d \cap[-r,r]^d$, it holds that for every $\varepsilon>0$,
\[
\lim_{k \to\infty} \mathbb{P} \bigl(\Lambda(\exp[k^{\Delta
^{-1}-\varepsilon}]) \subset\mathcal{B}_k \subset\Lambda(\exp
[k^{\Delta^{-1}+\varepsilon}])  \bigr) = 1.
\]
We note that an alternative proof of the statement
\[
\lim_{k \to\infty} \mathbb{P} \bigl(\mathcal{B}_k \subset\Lambda
(\exp[k^{\Delta^{-1}+\varepsilon}])  \bigr) = 1
\]
might be obtained by a slight change in the proof of Theorem \ref
{mainthm}(c): if we replace the definition ``$A_k(\varepsilon)$ is the
event that $\mathcal{B}_k$ contains a vertex at distance more than
$(1+ \varepsilon)^k$ from the origin'' by ``$A'_k(\varepsilon)$ is the
event that $\mathcal{B}_k$ contains a vertex at distance more than
$\exp[c k^{\gamma+\varepsilon}]$ from the origin,'' then the proof
essentially does not change.
%
%%\item For $\beta> 2d$, Coppersmith, Gamarnik and Sviridenko gave the
%following theorem \cite[Thm. 2.1(1)]{Copp02},
%%\begin{theorem}[Coppersmith, Gamarnik and Sviridenko, \cite{Copp02}]
%%Consider the supercritical homogeneous long-range percolation system
%as in \ref{nonexpothm}, with $\beta> 2d$. For $c_1> %%\frac{
%%\begin{equation}\label{corstate}
%%\lim_{\|x\| \to\infty} \mathbb{P} (D(0,x)>\|x\|^{1/c_1} ) =
%1.
%%\end{equation}
%%\end{theorem}
%%From Theorem \ref{nonexpothm2} we deduce the same statement as a
%corollary of Theorem \ref{coppthm} (Corollary \ref{newcorcop}), apart
%from the fact that we require $c_1> (\beta+1)/(\beta-2d)$. So
%Corollary \ref{newcorcop} is weaker than the result from
%$D(0,x)$ which is polynomially growing in $\|x\|$ is the same in both
%statements and Theorem \ref{nonexpothm2} also contains information
%about the rate at which $\mathbb{P} (D(0,x) \leq\|x\|^{1/c_1}
% )$ decreases to zero.
%%\item If a constant $\alpha>0$ exists such that $ \lim_{x
\eject
%s3 ###
\section[Proofs of Theorems 1.1 and 1.2]{Proofs of Theorems \protect\ref{mainthm}
and \protect\ref{nonexpothm}}
%s3.1 ###
\subsection[The $R_* = \infty$ regime: Proof of Theorem 1.1(a)]{The $R_* = \infty$ regime: Proof
of Theorem \protect\ref{mainthm}\textup{(a)}}

We consider long-range percolation with nonincreasing connection
function $\lambda(r) = r^{-\beta}L(r)$, where $L(r)$ is strictly
positive and slowly varying and $\beta< d$ or both $\beta= d$ and $\int
_{1}^{\infty} L(r)/ r \,dr = \infty$ hold.
We prove that in the cases under consideration in Theorem~\ref{mainthm}(a),
$ \sum_{x \in\mathbb{Z}^d \setminus\{0\}} p(0,x) =
\infty$ and, therefore,\vspace*{1pt} by the second Borel--Cantelli lemma (see,
e.g., \cite{Grim92}, page 288) we immediately obtain that $\mathbb
{P}(|\mathcal{B}_1| = \infty) =1$ a.s.

By \cite{Bing89}, Theorem 1.3.6, we know that for all $c>0$, $\lim_{r
\to\infty} r^{c}L(r) = \infty$ and so for both cases under
consideration and for all $R>0$, we obtain
\[
\int_R^{\infty} \lambda(r) r^{d-1} \,dr = \int_R^{\infty} (r^{d
-\beta} L(r))/r \,dr = \infty.
\]
Furthermore, note that for $x<1$, it holds that $1 - e^{-x} \geq x -
x^2/2 \geq x/2$
and that for large enough $r$, $\lambda(r)<1$ for both cases under
consideration.
Thus, constants $R>0$ and $c'>0$ exist such that
%
%e11 ###
%
\begin{eqnarray}
\sum_{x \in\mathbb{Z}^d \setminus0} p(0,x) &\geq& \sum_{x \in
\mathbb{Z}^d; \|x\|>R} p(0,x) \geq\frac{1}{2} \sum_{x \in
\mathbb{Z}^d; \|x\|>R} \lambda(\|x\|) \nonumber\\[-8pt]\\[-8pt]
&\geq& c' \int_{R+1}^{\infty}
r^{d-1} \lambda(r) \,dr= \infty,\nonumber
\end{eqnarray}
which proves that $|\mathcal{B}_1| = \infty$ a.s. in the regimes of
Theorem \ref{mainthm}(a).

%s3.2 ###
\subsection[The $R_* = 1$ regime: Proofs of Theorems 1.2 and 1.1(c)]{The $R_* = 1$ regime: Proofs
of Theorems \protect\ref{nonexpothm} and \protect\ref{mainthm}\textup{(c)}}

In this subsection, we prove that if $\lambda(r) = r^{-\beta}L(r)$,
with $\beta> d$ and $L(r)$ positive and slowly varying, then
$|\mathcal{B}_k|^{1/k} \to1$ a.s. and $\mathbb{E}(|\mathcal
{B}_k|^{1/k}|) \to1$ for $k \to\infty$.
To do this, we first show that Theorem \ref{nonexpothm} implies
Theorem \ref{mainthm}(c).
\begin{pf*}{Proof of Theorem \ref{mainthm}\textup{(c)}} Note that $|\mathcal
{B}_k| \geq1$ and, therefore, \mbox{$|\mathcal{B}_k|^{1/k} \geq1$} for all
$k$. So, $\liminf_{k \to\infty}(\mathbb{E}(|\mathcal{B}_k|))^{1/k}
\geq1$.
Furthermore, it is immediate from Theorem \ref{nonexpothm} that there
exists a constant $C$ such that
%
%e12 ###
%
\begin{eqnarray}
\mathbb{E}(|\mathcal{B}_k|) &=& \sum_{x \in V} \mathbb{P}\bigl(D(0,x) \leq k\bigr)
\nonumber\\[-8pt]\\[-8pt]
&\leq& \bigl(2K(k)+1\bigr)^d + \sum_{x \in V; \|x\| \geq K(k)} K(k)^{\beta'} \|x\|
^{-\beta'} \leq C K(k)^d.\nonumber
\end{eqnarray}
Because $\gamma<1$, we have $ \lim_{k \to\infty
}K(k)^{1/k} =1$. This implies that
\[
\limsup_{k \to
\infty}(\mathbb{E}(|\mathcal{B}_k|))^{1/k} = 1.
\]
Therefore, $R_* =1$.

From Theorem \ref{nonexpothm}, we obtain that if $d < \beta' < \beta
'' < \min(\beta,2d)$, then for all \mbox{$\varepsilon>0$}, there exists a
constant $N_1 = N_1(\varepsilon)$ such that for all $k >N_1$, it holds that $K(k) <
(1+\varepsilon)^k$.
Let $A_k = A_k(\varepsilon)$ be the event that $\mathcal{B}_k$ contains
a vertex at distance more than $(1+\varepsilon)^k$ from the origin.
For $k> N_1$, it holds that
%
%e13 ###
%
\begin{eqnarray}\qquad
\mathbb{P}(A_k) &\leq& \sum_{x \in\mathbb{Z}^d;\|x\|>(1+\varepsilon)^k}
[K(k)]^{\beta'}
\|x\|^{-\beta'} \leq c_1 \sum^{\infty}_{n=(1+\varepsilon)^k}
[K(k)]^{\beta'}n^{d-1-\beta'} \nonumber\\[-8pt]\\[-8pt]
&\leq& c_2[K(k)]^{\beta'} (1+\varepsilon)^{(d-\beta')k},\nonumber
\end{eqnarray}
where $c_1$ and $c_2$ are positive constants.
Note that for $d < \beta'<2d$, there exist constants $N_2>N_1$ and
$c_3>0$ such that for all $k>N_2$,
\[
c_2 \exp[c \beta' k^{\gamma}] (1+\varepsilon)^{(d-\beta')k} <
(1+\varepsilon)^{-c_3 k},
\]
by $\gamma:= \frac{\log(2d/\beta')}{\log(2)}<1$.
This implies that for every $\varepsilon>0$,
\[
\sum_{n=1}^{\infty} \mathbb{P}(A(n)) < N_2 + \sum_{n=N_2}^{\infty}
(1+\varepsilon)^{-c_3 n} < \infty
\]
and so $ \sum_{k=1}^{\infty} \mathbb{P}(|\mathcal
{B}_k|^{1/k}-1>\varepsilon) < \infty$,
which, in turn, implies, by the Borel--Cantelli lemma (see, e.g., \cite
{Grim92}, page 277), that Theorem \ref{mainthm}(c) holds.
\end{pf*}

Before providing the proof of Theorem \ref{nonexpothm}, we state a
useful lemma.
\begin{lemma}\label{hulplemma}
Consider the long-range percolation model defined in Section~\ref
{notation}, with vertex set $\mathbb{Z}^d$ and connection function
$\lambda(r)$ which satisfies $\lambda(r) < \alpha r^{-\beta''}$ for
all $r \geq1$ and a constant $\alpha>0$.  We
then have
\[
\mathbb{P}\bigl(D(0,x) \leq k\bigr) \leq\sum_{i=1}^{k} \mathbb{E}(|\mathcal
{B}_{i-1}|) \mathbb{E}(|\mathcal{B}_{k-i}|) \alpha(\|x\|/k)^{-\beta''}.
\]
\end{lemma}
\begin{pf}%\textbf{Proof of Lemma \ref{hulplemma}.}
If a self-avoiding path between $0$ and $x$ of length at most $k$
exists, then this path will contain at least one edge shared by
vertices at distance $\lceil\|x\|/k\rceil$ or more from each other.
Let $N(k,x)$ be the number of edges shared by vertices at distance at
least $\lceil\|x\|/k\rceil$ from each other that are contained in at
least one self-avoiding path between vertices $0$ and $x$ of length at
most $k$. For $1 \leq j \leq k$, let $N(k,x;j)$ be the number of edges
shared by vertices at distance at least $\lceil\|x\|/k\rceil$ from
each other that are contained as the $j$th edge in at least one
self-avoiding path from $0$ to $x$ of length at most $k$.

By Markov's inequality, we obtain
\[
\mathbb{P}\bigl(D(0,x) \leq k\bigr) \leq\mathbb{P}\bigl(N(k,x) \geq1\bigr) \leq\mathbb
{E}(N(k,x)) \leq\sum_{j=1}^{k} \mathbb{E}(N(k,x;j)).
\]
We further note that if we define $D(0,0) =0$, then, by observing that
for $x >0$, $p(x)< \lambda(x)$, it holds for $k \geq1$ that
\begin{eqnarray*}
&& \mathbb{E}(N(k,x;j))\\
&&\qquad \leq \mathop{\sum_{x_1,x_2 \in\mathbb{Z}^d}}_{\|x_1-x_2\|
\geq\lceil\|x\|/k \rceil} \mathbb{P}\bigl(D(0,x_1) = j-1\bigr) p(x_1,x_2)
\mathbb{P}\bigl(D(x_2,x) \leq k-j\bigr)\\
&&\qquad \leq \sum_{x_1,x_2 \in\mathbb{Z}^d} \mathbb{P}\bigl(D(0,x_1) =
j-1\bigr) p(\lceil\|x\|/k \rceil) \mathbb{P}\bigl(D(x_2,x) \leq k-j\bigr)\\
&&\qquad = \sum_{x_1 \in\mathbb{Z}^d} \mathbb{P}\bigl(D(0,x_1) = j-1\bigr)
p(\lceil\|x\|/k \rceil) \sum_{x_2 \in\mathbb{Z}^d}\mathbb
{P}\bigl(D(x_2,x) \leq k-j\bigr)\\
&&\qquad \leq \mathbb{E}(|\mathcal{B}_{j-1}|) \mathbb{E}(|\mathcal
{B}_{k-j}|) \alpha(\|x\|/k)^{-\beta''}.
\end{eqnarray*}
\upqed\end{pf}
\begin{pf*}{Proof of Theorem \ref{nonexpothm}}
We prove this theorem by using induction and Lemma \ref{hulplemma}.
Let the constant $c$ be such that
\[
c > \max \biggl(\frac{\log[\alpha]}{\beta'}, \frac{\beta''
+1}{(\beta''-\beta')\gamma} + \frac{3d \log[5] + 2 \log
[1+d/(\beta'-d)] + \log[\alpha]}{\beta''-\beta'} \biggr).
\]
Our induction hypothesis is that for all $j \leq k$ and $x \in\{x \in
\mathbb{Z}^d; x >K(j)\}$,
\[
\mathbb{P}\bigl(D(0,x) \leq j\bigr) \leq(K(j))^{\beta'} \|x\|^{-\beta'}.
\]
Note that the assumption holds for $k=0$ and $k=1$ because $K(0) =2$,
$K(1) = e^{c} +1$ and $c> (\beta')^{-1}\log[\alpha]$.
A straightforward computation yields that the induction hypothesis
implies that for all $j \leq k$,
\begin{eqnarray*}
\mathbb{E}(|\mathcal{B}_{j}|) & \leq& \bigl(2 K(j) +3\bigr)^d + \sum_{x \in
\mathbb{Z}^d; \|x\|>K(j)+1} (K(j))^{\beta'} \|x\|^{-\beta'}\\
& \leq& \bigl(2 K(j) +3\bigr)^d + \sum_{i= K(j) +2}^{\infty} 2d (2i
+1)^{d-1} (K(j))^{\beta'} i^{-\beta'}\\
& \leq& \bigl(2 K(j) +3\bigr)^d + 2d (K(j))^{\beta'} 3^{d-1} \int
_{K(j)}^{\infty} x^{d-\beta' - 1} \,dx\\
& = & \bigl(2 K(j) +3\bigr)^d + 2d (K(j))^{\beta'} 3^{d-1} (\beta' - d)^{-1}
{K(j)}^{d-\beta'}\\
& \leq&{K(j)}^d \bigl(5^d+ 2d 3^{d-1} (\beta' - d)^{-1}\bigr)\\
& \leq& 5^d\bigl(1 + d(\beta' - d)^{-1}\bigr) {K(j)}^d.
\end{eqnarray*}
We now observe that by Lemma \ref{hulplemma},
\begin{eqnarray*}
&& \mathbb{E}\bigl(N(k+1,x;j)\bigr)\\
&&\qquad \leq \mathbb{E}(|\mathcal{B}_{j-1}|) \mathbb{E}(|\mathcal
{B}_{k-j}|) \alpha\bigl(\|x\|/(k+1)\bigr)^{-\beta''}\\
&&\qquad \leq 5^{2d}  \biggl(1 + \frac{d}{\beta' - d} \biggr)^2
\bigl(K(j-1)K(k+1-j)\bigr)^d \alpha \biggl(\frac{k+1}{\|x\|} \biggr)^{\beta''}\\
&&\qquad \leq 5^{2d}  \biggl(1 + \frac{d}{\beta' - d} \biggr)^2 \bigl(e^{c
(j-1)^{\gamma}}+1\bigr)^d \bigl(e^{c (k+1-j)^{\gamma}}+1\bigr)^d \alpha \biggl(\frac
{k+1}{\|x\|} \biggr)^{\beta''}\\
&&\qquad \leq 5^{3d}  \biggl(1 + \frac{d}{\beta' - d} \biggr)^2
e^{2^{1-\gamma} d c k^{\gamma}} \alpha \biggl(\frac{k+1}{\|x\|}
\biggr)^{\beta''}\\
&&\qquad = 5^{3d}  \biggl(1 + \frac{d}{\beta' - d} \biggr)^2 (K(k))^{\beta
'} \alpha \biggl(\frac{k+1}{\|x\|} \biggr)^{\beta''},
\end{eqnarray*}
where we have used the fact that for $\gamma\leq1$ and $0 \leq y \leq
x$, $y^{\gamma} + (x-y)^{\gamma} \leq2 (x/2)^{\gamma}$, where the
right-hand side is equal to $x^{\gamma} \beta'/d$, by the definition
of $\gamma$.

Using this for $\|x\| >K(k+1)$, we obtain
\begin{eqnarray*}\qquad
&& \mathbb{P}\bigl(D(0,x) \leq k+1\bigr)\\
&&\qquad \leq \sum_{j=1}^{k+1} \mathbb{E}\bigl(N(k+1,x;j)\bigr)\\
&&\qquad \leq (k+1) 5^{3d}  \biggl(1 + \frac{d}{\beta' -
d} \biggr)^2 K(k)^{\beta'} \alpha \biggl(\frac{k+1}{\|x\|}
\biggr)^{\beta''}\\
&&\qquad = (k+1)^{\beta''+1} 5^{3d}  \biggl(1 + \frac
{d}{\beta' - d} \biggr)^2 K(k)^{\beta'} \alpha\|x\|^{-(\beta''-\beta
')} \|x\|^{-\beta'}\\
&&\qquad \leq (k+1)^{\beta''+1} 5^{3d}  \biggl(1 + \frac
{d}{\beta' - d} \biggr)^2
\alpha\bigl(K(k+1)\bigr)^{-(\beta''-\beta')}\\
&&\qquad\quad{}\times
\bigl(K(k+1)\bigr)^{\beta'} \|x\|^{-\beta'}.
\end{eqnarray*}
Because
\[
c > \frac{\beta'' +1}{(\beta''-\beta')\gamma} + \frac{3d \log[5]
+ 2 \log[1+d/(\beta'-d)] + \log[\alpha]}{\beta''-\beta'}
\]
and $\mathbb{P}(D(0,x)\leq k+1) \leq1$, we deduce, after some
straightforward computations,
that
\[
(k+1)^{\beta''+1} 5^{3d}\bigl(1 + d(\beta' - d)^{-1}\bigr)^2 \alpha\exp
[-(\beta''-\beta') c (k+1)^{\gamma}] \leq1
\]
for $k \geq0$ and that if the induction hypothesis holds, then it also
holds that if $\|x\|> K(k+1) $, then
\[
\mathbb{P}\bigl(D(0,x) \leq k+1\bigr) \leq\bigl(K(k+1)\bigr)^{\beta'} \|x\|^{-\beta'}.
\]
This proves the theorem.
\end{pf*}

Theorem \ref{nonexpothm} contains a part of Theorem \ref{biskupthm}
as the following corollary.
\begin{corollary}\label{newcorbisk}
Consider homogeneous long-range percolation on $\mathbb{Z}^d$, as in
Section \ref{nonexpothm}, with nonincreasing connection function
$\lambda(r) = r^{-\beta}L(r)$, where $d< \beta< 2d$ and $L(r)$ is
positive and slowly varying. For $\Delta= \frac{\log[2]}{\log
[2d/\beta]}$ and every $\varepsilon>0$, we have
%
%e14 ###
%
\begin{equation}\label{corstate2}
\lim_{\|x\| \to\infty} \mathbb{P} \biggl(\Delta-\varepsilon< \frac
{\log[D(0,x])}{\log[\log[\|x\|]]} \biggr) = 1.
\end{equation}
\end{corollary}
\begin{pf}
Observe that (\ref{corstate2}) can be rewritten as
\[
\lim_{\|x\| \to\infty} \mathbb{P} \bigl(D(0,x) \leq(\log[\|x\|
])^{\Delta-\varepsilon} \bigr) = 0.
\]
We choose $\beta'< \beta$ such that
\[
(\Delta- \varepsilon) \gamma=  \biggl(\frac{\log[2]}{\log[2d/\beta]}
- \varepsilon \biggr) \biggl(\frac{\log[2]}{\log[2d/\beta']}
\biggr)^{-1} < 1,
\]
which can be done for every $\varepsilon> 0$.
By substituting $k = (\log[\|x\|])^{\Delta- \varepsilon}$ into
\[
\mathbb{P}\bigl(D(0,x) \leq k\bigr) \leq[K(k)]^{\beta'} \|x\|^{-\beta'}
\qquad\mbox{for $x \in\{x \in\mathbb{Z}^d; \|x\| > K(k)\} $},
\]
and $K(k) = 1 + \exp[c k^{\gamma}]$, we obtain that
\[
\mathbb{P}\bigl(D(0,x) \leq(\log[\|x\|])^{\Delta- \varepsilon}\bigr) \leq
\bigl(1+\exp\bigl[c(\log[\|x\|])^{(\Delta-\varepsilon)\gamma}\bigr]\bigr)^{\beta'} \|x\|
^{-\beta'}
\]
for $x \in\{x \in\mathbb{Z}^d; \|x\| > 1+\exp[c (\log[\|x\|
])^{(\Delta-\varepsilon)\gamma}]\}$.
If $\|x\|$ is large enough, then
\[
\|x\|^{1/2} > 1+\exp\bigl[c(\log[\|x\|])^{(\Delta-\varepsilon)\gamma}\bigr]
\quad\Leftrightarrow\quad\|x\|^{1/2} > 1+\|x\|^{c(\log[\|x\|])^{(\Delta-\varepsilon
)\gamma-1}}
\]
holds. Therefore, for $\|x\| \to\infty$,
\begin{eqnarray*}
\mathbb{P}\bigl(D(0,x) \leq(\log[\|x\|])^{\Delta- \varepsilon}\bigr) & \leq&
\bigl(1+\exp\bigl[c\log[\|x\|]^{(\Delta-\varepsilon)\gamma}\bigr]\bigr)^{\beta'} \|x\|
^{-\beta'}\\
& = & \bigl(\|x\|^{-1}+\|x\|^{c(\log[\|x\|])^{(\Delta-\varepsilon)\gamma
-1}-1}\bigr)^{\beta'} \\
& \to& 0,
\end{eqnarray*}
which proves the corollary.
\end{pf}

%s3.3 ###
\subsection[The nontrivial $R_*$ regime: Proof of Theorem 1.1(b)]{The nontrivial $R_*$ regime: Proof
of Theorem \protect\ref{mainthm}\textup{(b)}}

In this subsection, we consider percolating homogeneous long-range
percolation with nonincreasing connection function $\lambda(r) =
r^{-d}L(r)$, where $L(r)$ is nonnegative, slowly varying and satisfies
\[
\int_{1}^{\infty} r^{-1}L(r)\,dr < \infty \quad\mbox{and}\quad
-\int_{R}^{\infty} \frac{\log[L(r)]}{r (\log[r])^2} \,dr< \infty
\]
for some constant $R>0$.
We investigate the growth behavior of $|\mathcal{B}_k|$ for \mbox{$k \to
\infty$}. In particular, we show that $\lambda(r)$ exists, satisfying
these conditions, such that ${\lim_{k \to\infty}}|\mathcal
{B}_k|^{1/k} >1$ with positive probability.

In the first subsection, we provide a straightforward and almost
trivial proof for the upper bounds of the growth of $|\mathcal{B}_k|$
given in Theorem \ref{mainthm}(b). After that, we provide some useful
lemmas that will be used in the proof of the lower bound of the growth
of $|\mathcal{B}_k|$. We then give an outline of the proof and, in the
final subsection, the full proof of $\lim_{k \to\infty} \mathbb{P}(
|\mathcal{B}_k|^{1/k} > a_1 | 0 \in\mathcal{C}_{\infty}) =1$ for
some $a_1>1$ is given. In this proof, renormalization arguments are used.

%s3.3.1 ###
\subsubsection{The upper bound: Proof of $\overline{R}_* < \infty$}
\begin{lemma}\label{brwupperlem}
Consider the percolating homogeneous long-range percolation model, as
defined in Section \ref{notation}, with vertex set $\mathbb{Z}^d$ and
nonincreasing connection function $\lambda(r) = r^{-d} L(r)$, where
$L(r)$ is positive, slowly varying and satisfies condition (\ref
{sumable}). There exists a positive and finite constant $a_2$ such that
$\overline{R}_* < a_2$ and
\[
\lim_{k \to\infty} \mathbb{P}( |\mathcal{B}_k|^{1/k} < a_2) =1.
\]
\end{lemma}
\begin{pf}
Assign an independent Poisson process to each pair of vertices in
$\mathbb{Z}^d$, representing the contacts between the pair of
vertices. The density of the Poisson process of vertices at distance $r
> 0$ is $\lambda(r)$.
We observe that the probability that at least one contact is made
between two vertices at distance $r$ in the interval $(0,1)$ is~$p(r)$.
If the pairs of individuals that make at least one contact in the
interval $(0,1)$ are joined by an edge, then the long-range percolation
graph under consideration is re-obtained.

We obtain, after some basic computations, that
\[
 \sum_{x \in\mathbb{Z}^d \setminus\{y\}} \lambda(\|
x-y\|)< \infty
\]
for all $y \in\mathbb{Z}^d$.
%%\begin{eqnarray*}
%%&  & \sum_{v \in\mathbb{Z}^d; |v-w|>K_1} \lambda(|w-v|)\\
%%\ & \leq& \sum_{k=K_1+1}^{\infty} 2d (2k+1)^{d-1} \lambda(k)\\
%%\ & \leq& 2d[(2K_1 +3)^{d-1} \lambda(\lfloor K_1 \rfloor+1) +
%(2K_1+5)^{d-1} \lambda(\lfloor K_1 \rfloor+2)] + \sum_{k=K_1+3}^{
%%\ & \leq& 2d[(2K_1 +3)^{d-1}\lambda(\lfloor K_1 \rfloor+1) +
%(2K_1+5)^{d-1} \lambda(\lfloor K_1 \rfloor+2)] + \sum_{k=K_1+3}^{
%%\ & \leq& 2d[(2K_1 +3)^{d-1}\lambda(\lfloor K_1 \rfloor+1) +
%(2K_1+5)^{d-1} \lambda(\lfloor K_1 \rfloor+2)] + 2d 3^{d-1}
%%\ & \leq& 2d[(2K_1 +3)^{d-1}\lambda(\lfloor K_1 \rfloor+1) +
%(2K_1+5)^{d-1} \lambda(\lfloor K_1 \rfloor+2)] + 2d 6^{d-1}
%%\ & \leq& \infty.
%%\end{eqnarray*}
It is straightforward to couple the $k$-ball, $|\mathcal{B}_k|$, to
the number of individuals in the first $k$ generations of a
supercritical branching random walk with a Poisson-distributed
offspring size distribution and $\overline{R}_* < a_2$ for some
$a_2>1$ follows immediately. The branching random walk is the process
in which, initially, one individual (or particle) lives at the origin.
This individual stays there forever, although it can only give birth to
new individuals during the first time unit of its life. This individual
gives birth to individuals at vertex $x$ according to a Poisson process
with rate $\lambda(\|x\|)$.

The set $\mathcal{B}_k$ is created by killing upon birth all
individuals that are born on a vertex that is already occupied by
another individual. From the theory of branching processes \cite
{Jage75}, we know that there exist a random variable, $W$, which is
almost surely finite, and a constant $a'_2$ such that for all $k \in
\mathbb{N}_{+}$, the number of individuals in the first $k$
generations of such a branching random walk is a.s. bounded above by $W
(a'_2)^k$.

In the coupled process, $|\mathcal{B}_k|$ is bounded above by the
number of individuals in the first $k$ generations of the branching
random walk,
which proves that $ \lim_{k \to\infty} \mathbb
{P}(|\mathcal{B}_k|^{1/k} < a_2) =1$ for $a_2>a'_2$.
\end{pf}

%s3.3.2 ###
\subsubsection{The lower bound, preliminary lemmas and definitions}

In order to prove that for the given connection function,
$ \lim_{k \to\infty} \mathbb{P}(a_1 < |\mathcal
{B}_k|^{1/k}| 0 \in\mathcal{C}_{\infty}) =1$ holds, we need the
following lemmas.
\begin{lemma}\label{assumpequi}
For positive, slowly varying $L(x)$ which is nonincreasing and strictly
less than 1 on $[K,\infty)$ for some $K>1$, condition (\ref
{annoyass}) is equivalent to the following condition:
for every $\delta> 0$, some $K_1 := K_1(\delta) > K$ exists such that
for $r>K_1$,
\[
\sum_{k=1}^{\infty} \frac{\log[L(r^{2^k})]}{2^{k}(\log[r])}
>-\delta
\]
and, therefore,
\[
\prod_{k=1}^{\infty}[L(r^{2^k})]^{2^{-k}}>r^{-\delta}.
\]
\end{lemma}
\begin{lemma}\label{positg}
For positive, slowly varying $L(x)$ which is nonincreasing on
$[K,\infty)$ for some $K>1$ and satisfies condition (\ref{sumable}),
there exists $K_2 \geq K$ such that $L(r)<1$ for all $r \geq K_2$.
\end{lemma}
\begin{lemma}\label{nontrivpcor}
Consider a percolating homogeneous long-range percolation graph, as in
Section \ref{notation}, with nonincreasing connection function
$\lambda(r) = r^{-d}L(r)$, where $L(r)$ is positive and slowly varying.
Let $\mathcal{C}'_K$ be the largest cluster of $G$ restricted to
vertices in $V_K := \mathbb{Z}^d \cap[\lceil-K/2 \rceil,\lceil K/2
\rceil)^d$ and edges shared by vertices both in $V_K$.
For every $\varepsilon>0$, there exist numbers $\rho>0$ and $K_3
=K_3(\varepsilon) < \infty$ such that for every $r \geq K_3$,
%
%e16 ###
%e15 ###
%
\begin{eqnarray}
\mathbb{P}(|\mathcal{C}'_r| < \rho r^d) & \leq& \varepsilon,\\
\mathbb{P}(|\mathcal{C}'_r| < \rho r^d, 0 \in\mathcal{C}'_r| 0 \in
\mathcal{C}_{\infty}) & \leq& \varepsilon.
\end{eqnarray}
\end{lemma}
\begin{pf*}{Proof of Lemma \ref{assumpequi}}
Since $L(x)$ is nonincreasing and strictly less than~1 for $x \geq K$,
it is\vspace*{-4pt} enough to show that $- \int_K^{\infty} \frac{\log[L(x)]}{x
(\log[x])^2} \,dx < \infty$ is equivalent to $-\sum_{k=1}^{\infty}
\frac{\log[L(K^{2^k})]}{2^{k}(\log[K])} < \infty$.\vspace*{2pt}

From \cite{Bing89}, we know that for slowly varying, eventually
decreasing $L(x)$, there exist a function $\delta(x)$, converging to a
finite number, and a nonnegative function~$\varepsilon(x)$, converging to
$0$ for $x \to\infty$, such that
\[
L(x) = \exp\biggl[\delta(x) - \int_K^{x} \frac{\varepsilon(t)}{t} \,dt\biggr].
\]
Substituting this in gives
%
%e17 ###
%
\begin{eqnarray}
&&\int_K^{\infty} \frac{\log[L(x)]}{x(\log[x])^2} \,dx \nonumber\\
&&\qquad= \int
_K^{\infty} \frac{\delta(x) - \int_K^{x} {\varepsilon(t)}/{t}
\,dt}{x (\log[x])^2} \,dx \nonumber\\[-8pt]\\[-8pt]
&&\qquad= \int_K^{\infty} \frac{\delta(x)}{x (\log[x])^2} \,dx - \int
_K^{\infty} \int_t^{\infty} \frac{\varepsilon(t)}{t} \frac{1}{x
(\log[x])^2} \,dx \,dt \nonumber\\
&&\qquad= \int_K^{\infty} \frac{\delta(x)}{x (\log[x])^2} \,dx - \int
_K^{\infty} \frac{\varepsilon(t)}{t \log[t]} \,dt.\nonumber
\end{eqnarray}
Since $\delta(x)$ converges to a finite number and $0 < L(x) < 1$ for
$x \geq K$, it follows that $\delta(x)$ is bounded away from infinity
and the first term is finite.
Thus, we obtain that $- \int_K^{\infty} \frac{\log[L(x)]}{x (\log
[x])^2} \,dx < \infty$ is equivalent to $\int_K^{\infty} \frac
{\varepsilon(t)}{t \log[t]} \,dt< \infty$.

Note that
\[
\sum_{k=1}^{\infty} \frac{\log[L(K^{2^k})]}{2^{k}(\log[K])} = \sum
_{k=1}^{\infty} \frac{\delta(K^{2^k})}{2^{k}(\log[K])} - \sum
_{k=1}^{\infty} \frac{\int_K^{K^{2^k}} {\varepsilon(t)}/{t}
\,dt}{2^{k}(\log[K])}
\]
and
%
%e18 ###
%
\begin{eqnarray}\label{snel}
\sum_{k=1}^{\infty} \frac{\int_K^{K^{2^k}} {\varepsilon(t)}/{t}
\,dt}{2^{k}(\log[K])} &=& \sum_{k=1}^{\infty} \sum_{l=1}^k
\frac{\int_{K^{2^{l-1}}}^{K^{2^l}} {\varepsilon(t)}/{t}
\,dt}{2^{k}(\log[K])}\nonumber\\
&=& \sum_{l=1}^{\infty} \sum_{k=l}^{\infty}
\frac{\int_{K^{2^{l-1}}}^{K^{2^l}} {\varepsilon(t)}/{t}
\,dt}{2^{k}(\log[K])} \\
&=& \sum_{l=1}^{\infty}
\frac{\int_{K^{2^{l-1}}}^{K^{2^l}} {\varepsilon(t)}/{t}
\,dt}{2^{l-1}(\log[K])}.\nonumber
\end{eqnarray}
The final term of (\ref{snel}) is bounded below by
\[
\sum_{l=1}^{\infty} \int_{K^{2^{l-1}}}^{K^{2^l}} \frac{\varepsilon
(t)}{t \log[t]} \,dt = \int_{K}^{\infty} \frac{\varepsilon(t)}{t \log
[t]} \,dt.
\]
Similarly,\vspace*{-4pt} we deduce that this term is bounded above by $2 \int
_{K}^{\infty} \frac{\varepsilon(t)}{t \log[t]} \,dt$. Furthermore, $\sum
_{k=1}^{\infty} \frac{\delta(K^{2^k})}{2^{k}(\log[K])}$ is
finite\vspace*{1pt}
since $\delta(x)$ converges and is bounded away from infinity. Therefore,
\[
-\sum_{k=1}^{\infty} \frac{\log[L(K^{2^k})]}{2^{k}(\log[K])} <
\infty\quad\Leftrightarrow\quad\int_K^{\infty} \frac{\varepsilon(t)}{t \log
[t]} \,dt< \infty
\]
and the lemma follows.
\end{pf*}
\begin{pf*}{Proof of Lemma \ref{positg}}
Because $L(x)$ is not increasing for $x>K$, we know that $\lim_{x \to
\infty} L(x)$ exists. If $\lim_{x \to\infty} L(x)>0$, then
$\int_{K}^{\infty} L(x) x^{-1} \,dx = \infty$ and this violates
condition (\ref{sumable}). Together with the assumption that $L(x) >0$
for all \mbox{$x>0$}, this leads to $\lim_{x \to\infty} L(x)= 0$ and
therefore there exists $K_2 \in\mathbb{R}_+$ such that $L(x)<1$ for
all $x>K_2$.
\end{pf*}
\begin{pf*}{Proof of Lemma \ref{nontrivpcor}} This follows immediately
from Theorem 3.2 and Corollary 3.3 in \cite{Bisk04}, together with the
fact that $\mathbb{P}(0 \in\mathcal{C}_{\infty}) > 0$.
\end{pf*}

In the remaining proof of Theorem \ref{mainthm}(b), we use a
construction for which the following definitions are needed.

We define a hierarchy of blocks of vertices in $\mathbb{Z}^d$ as
follows. For every $i \in\mathbb{N}$ and every $\bar{n} \in\mathbb{Z}^d$,
we define $\Lambda_i(\bar{n}) := \mathbb{Z}^d \cap
[-(l_i-1)/2,(l_i-1)/2] + \bar{n} l_i$,
where $l_i := (l_0)^{2^i}$, $l_0$ is odd and
\[
l_0 > \max\biggl[\biggl(\frac{100}{\rho^2}\biggr)^{1/(d-2/5)}, \frac{100}{16d},
K_1(1/5), K_2,K_3(1/25)\biggr].
\]
The constants $\rho,K_1(1/5),K_2$ and $K_3(1/25)$ are as in the
preceding lemmas.
We say that $\Lambda_i(\bar{n})$ is a \textit{level} $i$ \textit{block} and
note that for every $i \in\mathbb{N}$, the level $i$ blocks form a
partition of $\mathbb{Z}^d$. Every level $i$ block is entirely
contained in a level $i+1$ block and every level $i+1$ block contains
$(l_0)^{2^i d}$ level $i$ blocks. We use $\bar{n}_i(x)$ to denote the
index of the level $i$ block containing vertex $x$, that is, for $x \in
V$, we define $\bar{n}_i(x) \in\mathbb{Z}^d$ such that $x \in
\Lambda_i(\bar{n}_i(x))$.

Let $G$ be the long-range percolation graph under consideration and let
$G_i(\bar{n})$ be defined as $G$ restricted to $\Lambda_i(\bar{n})$,
that is, $G_i(\bar{n})$ is the graph consisting of vertex set $\Lambda
_i(\bar{n})$ and those edges of $G$ for which both end-vertices are in
$\Lambda_i(\bar{n})$.
Let $\mathcal{D}_i(x)$ be the set of vertices in $\Lambda_i(\bar
{n}_i(x))$ that are within graph distance
\[
h_i :=(l_0)^d 2^i + {2^i-1} = \bigl((l_0)^d+1\bigr)2^i-1
\]
of $x$ in the graph $G_i(\bar{n}_i(x))$.

A vertex $x \in\mathbb{Z}^d$ is said to be \textit{good up to level
$0$} if
\[
|\mathcal{D}_0(x)| \geq m_0 := \rho(l_0)^d.
\]

For $x \in\mathbb{Z}^d$ and $S \subset\mathbb{Z}^d$, let $x
\leftrightarrow S$ denote the event that there is a vertex $y \in S$
such that $\langle x,y\rangle  \in E$. Furthermore, let
\[
\bar{\mathcal{D}}_{i+1}(x):= \bigl\{y \in\mathbb{Z}^d \setminus\bigl(\Lambda
_i(0) \cup\Lambda_i(\bar{n}_i(x))\bigr)| y \leftrightarrow\mathcal
{D}_i(x),\mbox{ $y$ good up to level $i$}\bigr\}
\]
be the set of vertices not in $(\Lambda_i(0) \cup\Lambda_i(\bar
{n}_i(x)))$ that share an edge with vertices in $\mathcal{D}_i(x)$ and
that are good up to level $i$.
A vertex $x \in\mathbb{Z}^d$ is \textit{good up to level $i+1$} if
$x$ is good up to level $i$ and if
%
%e19 ###
%
\begin{eqnarray}
A_{i+1}(x)
:\!&=& |\{\bar{n} \in\mathbb{Z}^d|\Lambda_i(\bar{n}) \subset\Lambda
_{i+1}(\bar{n}_{i+1}(x)),\nonumber\\
&&\hspace*{20.4pt}\exists y \in\Lambda_i(\bar{n}),\mbox{ s.t. $y \in\bar{\mathcal
{D}}_{i+1}(x)$}\}| \\
&\geq& m_{i+1} := c_0 L(l_{i+1}) M_i.\nonumber
\end{eqnarray}
Here,\vspace*{1pt} $c_0 :=2 \rho/25$, $M_i := \prod_{j=0}^i m_j$ and for $i \in
\mathbb{N}_+$, the constants $m_i$ are defined recursively.
In words, this means that the number of level $i$ blocks in $\Lambda
_{i+1}(\bar{n}_{i+1}(x)\setminus(\Lambda_i(0) \cup\Lambda_i(\bar
{n}_i(x)))$ that contain at least one vertex that shares an edge with a
vertex in $D_i(x)$ is at least $m_{i+1}$.
Some algebra gives that $m_{i+1} = L(l_{i+1}) (c_0 m_0 \prod_{j=1}^i
[L(l_j)]^{2^{-j}})^{2^i} = L(l_{i+1}) (c_0 \rho(l_0)^d \prod_{j=1}^i
[L(l_j)]^{2^{-j}})^{2^i}$ for\vspace*{1pt} $i \in\mathbb{N}$. Since $l_0>
K_1(1/5)$ and $l_0>K_2$, Lemmas \ref{assumpequi} and \ref{positg}
give that
\[
\prod_{j=1}^i [L(l_j)]^{2^{-j}} \geq\prod_{j=1}^{\infty}
[L(l_j)]^{2^{-j}} \geq(l_0)^{-1/5},
\]
while Lemma \ref{positg} gives
\[
\prod_{j=1}^i [L(l_j)]^{2^{-j}} \leq1.
\]
Combining these observations gives that for $i \in\mathbb{N}_+$,
\[
L((l_0)^{2^i}) (c_0 \rho(l_0)^{d-1/5})^{2^{i-1}} \leq m_i \leq
L((l_0)^{2^i}) (c_0 \rho(l_0)^{d})^{2^{i-1}}.
\]
Since $L(x)$ is slowly varying, this implies that there exist constants
$1 < c'_0 < c''_0 < \infty$ and $0< \hat{c}'_0< \hat{c}''_0<\infty$
such that $\hat{c}'_0(c'_0)^{2^i} < m_i < \hat{c}''_0(c''_0)^{2^i}$
holds for all $i \in\mathbb{N}_+$.

A vertex is \textit{ultimately good} if it is good up to every level
$i \in\mathbb{N}$.

%s3.3.3 ###
\subsubsection[Outline of proof of lower bound in Theorem 1.1(b)]{Outline of proof of
lower bound in Theorem \protect\ref{mainthm}\textup{(b)}}
As may be guessed from the definitions above, the proof will follow a
renormalization scheme. The following steps are taken.
\begin{itemize}
\item We observe that if $x$ is good up to level $i+1$, then $|\mathcal
{D}_{i+1}(x)| \geq m_{i+1}|\mathcal{D}_{i}(x)|$ and $x$ is good up to
all levels $0 \leq j \leq i$. Therefore,
%
%e20 ###
%
\begin{eqnarray}\label{kolgood}
|\mathcal{D}_{i}(x)| &\geq& M_i = m_{i+1} (c_0 L(l_{i+1})^{-1})
\nonumber\\
&=& \frac{1}{c_0}  \Biggl(c_0 \prod_{j=1}^{i} [L(l_j)]^{2^{-j}} m_0
 \Biggr)^{2^i}
= \frac{1}{c_0}  \Biggl(c_0 \rho\prod_{j=1}^{i} [L(l_j)]^{2^{-j}}
(l_0)^d  \Biggr)^{2^i} \\
&\geq&\frac{1}{c_0}  \Biggl(c_0 \rho\prod_{j=1}^{\infty}
[L(l_j)]^{2^{-j}} (l_0)^d  \Biggr)^{2^i}
\geq\frac{1}{c_0}  (c_0 \rho(l_0)^{d-1/5}  )^{2^i}.\nonumber
\end{eqnarray}
Here, we have used the fact that $l_0 > K_2$. Note that, by $l_0 >
(100/\rho^2)^{1/(d-2/5)}$, we have
%
%e21 ###
%
\begin{equation}\label{usefulbound}
c_0 \rho(l_0)^{d-1/5} > (2 \rho^2/25) (100/\rho^2)^{(d-1/5)/(d-2/5)}
> 8.
\end{equation}
\item Recall that $\mathcal{B}_{k}(x)$ is the set of vertices in
$\mathbb{Z}^d$ within graph distance $k$ (in $G$) of $x$. Note that
$\mathcal{D}_i(x) \subset\mathcal{B}_{h_i}(x)$.
We show that if $x$ is ultimately good, then $|\mathcal
{B}_{h_i}(x)|^{1/{h_{i+1}}} > a'_1$ for some $a'_1>1$, which, in turn,
implies that if $x$ is ultimately good, then $|\mathcal
{B}_{k}(x)|^{1/k} > a_1$ for all $k \geq1$ and some $a_1>1$.
\item We show\vspace*{1pt} that $l_0$ is large enough to guarantee that the
probability that $x$ is ultimately good is positive and $\rho c_0
(l_0)^{d-1/5}>1$.
\item We use a zero--one law to prove that the number of ultimately good
vertices is infinite.
\item Finally, we show that $|\mathcal{B}_{j}|^{1/j} := |\mathcal
{B}_{j}(0)|^{1/j} > a_1$ if $0 \in\mathcal{C}_{\infty}$.
\end{itemize}

%s3.3.4 ###
\subsubsection{Proof of $\underline{R}{}_*>1$}

We are now ready to state a lemma which will lead to the proof of
Theorem \ref{mainthm}(b).
\begin{lemma}\label{kolmouse}
Consider a\vspace*{1pt} percolating homogeneous long-range percolation model, as
defined in Section \ref{notation}, with vertex set $\mathbb{Z}^d$ and
connection function $\lambda(r) = r^{-d}L(r)$, where $L(r)$ is
positive, slowly varying, decreasing on $[K,\infty)$ for some $K>0$
and satisfies (\ref{sumable}) and (\ref{annoyass}). If $\rho,c_0$
and $l_0$ are as above, then the number of ultimately good vertices in
$\mathbb{Z}^d$ is a.s. infinite.
\end{lemma}
\begin{pf*}{Proof of Theorem \ref{mainthm}\textup{(b)}}
Note that $\rho c_0(l_0)^{d-1/5}>1$, by (\ref{usefulbound}), and if
$x$ is ultimately good, then (\ref{kolgood}) implies that for
$2^i (1+(l_0)^{d}) \leq k < 2^{i+1} (1+(l_0)^{d})$,
\[
|\mathcal{B}_k(x)| \geq\frac{1}{c_0}  (\rho
c_0(l_0)^{d-1/5} )^{2^i} \geq c_0^{-1} (\rho
c_0(l_0)^{d-1/5})^{(1 + (l_0)^d)^{-1}k/2},
\]
where we have used the fact that $l_0 > \max(K_1(1/5),K_2)$.

Lemma \ref{kolmouse} implies that there is at least one ultimately
good vertex in $\mathbb{Z}^d$.
By the construction of $\mathcal{D}_i(x)$, it is clear that an
ultimately good vertex is in an infinite cluster of $G$. By the
uniqueness of the infinite cluster of $G$, we know that conditioned on
$\{0 \in\mathcal C_{\infty}\}$, the random variable $Y := \min\{
D(0,x); x \in\mathbb{Z}^d, \mbox{$x$ is ultimately good}\}$ is
a.s. finite. Therefore,
\begin{eqnarray*}
(\mathcal{B}_{k+Y})^{1/(k+Y)} & \geq&  \bigl(((c_0)^{-1} (\rho
c_0(l_0)^{d-1/5}))^{(1 + (l_0)^d)^{-1}k/2} \bigr)^{1/(k+Y)} \\
& \geq& \bigl((c_0)^{-1/(k+Y)} (\rho c_0(l_0)^{d-1/5})\bigr)^{(1 +
(l_0)^d)^{-1}k/(2 (k+ Y))},
\end{eqnarray*}
which converges to $a'_1 := (\rho c_0
(l_0)^{d-1/5})^{(2+2(l_0)^{d})^{-1}} >1$ and, therefore, there exists a
constant $a_1>1$ such that
\[
\lim_{k \to\infty} \mathbb{P}(a_1 < |\mathcal{B}_k|^{1/k}| 0 \in
\mathcal{C}_{\infty}) =1,
\]
which proves the theorem.
\end{pf*}

For the proof of Lemma \ref{kolmouse}, we need a bound for
\[
\mathbb{P}(\mbox{$x$ is good up to level $i+1| x$ is good up to level $0$}).
\]
We obtain this bound by using the following lemma.
\begin{lemma}\label{hierarchlem1}
Consider a percolating homogeneous long-range percolation model, as
defined in Section \ref{notation}, with vertex set $\mathbb{Z}^d$ and
connection function $\lambda(r) = r^{-d}L(r)$, where $L(r)$ is
positive, slowly varying, nonincreasing on $[K,\infty)$ for some $K>0$
and satisfies (\ref{sumable}) and (\ref{annoyass}). If $\rho,c_0$
and $l_0$ are as above, then, for $i \in\mathbb{N}$,
\[
\mathbb{P}(\mbox{$x$ is good up to level $i+1| x$ is good up to level
$i$}) \geq1-4^{-2^i}.
\]
\end{lemma}
\begin{pf}
If we assume that the
statement holds for $j<i$, then
%
%e22 ###
%
\begin{eqnarray}\label{lemmacomp}\quad
&&\mathbb{P}(\mbox{$x$ is good up to level $i | x$ is good up to level
$0$}) \nonumber\\
&&\qquad\geq1- \sum_{j=0}^{i-1} \mathbb{P}(\mbox{$x$ is not good up to
level $j+1 | x$ is good up to level $j$}) \\
&&\qquad\geq1- \sum_{j=0}^{i-1} 4^{-2^j} \geq1- \sum_{j=0}^{i-1}
4^{-(j+1)} = 1 - \frac{1-4^{-i}}{3} \geq2/3.\nonumber
\end{eqnarray}
Furthermore, note that if the random variable $X$ is binomially
distributed with parameters $n$ and $p$, then, by Chebyshev's
inequality, we have
%
%e23 ###
%
\begin{equation}\label{binafs}
\mathbb{P} \biggl(X <\frac{\mathbb{E}(X)}{2} \biggr) \leq\frac{4
\operatorname{Var}(X)}{(\mathbb{E}(X))^2} = \frac{4 (1-p)}{np} \leq\frac{4}{np}.
\end{equation}
Observe that if $x$ and $y$ are not in the same level $i$ block, then
the events $\{\mbox{$y$ is good up to level $i$}\}$ and $\{y
\leftrightarrow\mathcal{D}_i(x)\}$ are independent, because different
edges are involved.
We already know, by (\ref{kolgood}), that if $x$ is good up to level
$i$, then
%
%e24 ###
%
\begin{equation}\label{binafs2}
|\mathcal{D}_{i}(x)| \geq\frac{1}{c_0}  \Biggl(c_0 \rho\prod
_{j=1}^{i} [L(l_j)]^{2^{-j}} (l_0)^d  \Biggr)^{2^i}.
\end{equation}
Furthermore, all vertices in $\Lambda_{i+1}(\bar{n}_{i+1}(x))
\setminus(\Lambda_i(0) \cup\Lambda_{i}(\bar{n}_{i}(x)))$ have
probability exceeding $1 - \exp[- |\mathcal{D}_{i}(x)| \lambda
(l_{i+1})]$ to share an edge with a vertex in $\mathcal{D}_{i}(x)$.
Therefore, the probability that a given level $i$ block,
\[
\Lambda_{i}(\bar{n}') \subset\Lambda_{i+1}(\bar{n}_{i+1}(x)) \setminus
\bigl(\Lambda_i(0) \cup\Lambda_{i}(\bar{n}_{i}(x))\bigr),
\]
contains a vertex (say $y$) that is good up to level $i$ and shares an
edge with a vertex in $\mathcal{D}_{i}(x)$ is bounded below by
\begin{eqnarray*}
&& \mathbb{P} \bigl(\mbox{$y$ is good up to level $0$}| \mbox
{$y$ is chosen uniformly at random from $\Lambda_{i}(\bar
{n}_{i}(y))$} \bigr) \\
&&\qquad{} \times \mathbb{P} (\mbox{$y$ is good up to level
$i$}|\mbox{$y$ is good up to level $0$} ) \\
&&\qquad{} \times \mathbb{P}^* \Biggl(\Lambda_{i}(\bar{n}_{i}(y))
\leftrightarrow\mathcal{D}_i(x)||\mathcal{D}_{i}(x)| = \frac
{1}{c_0}  \Biggl(c_0 \rho\prod_{j=1}^{i} [L(l_j)]^{2^{-j}} (l_0)^d
 \Biggr)^{2^i}  \Biggr).
\end{eqnarray*}
Here, $\mathbb{P}^*$ is the product measure for which a pair of
vertices $x,y \in\mathbb{Z}^d$ share an edge with probability
$1-e^{-\lambda(l_{R(x,y)})}$, where
\[
R(x,y) = \inf\{i \in\mathbb{N}; y \in\Lambda_i(\bar{n}_i(x))\}.
\]
Note that
\begin{eqnarray*}
& & \mathbb{P}^* \Biggl(\Lambda_{i}(\bar{n}_{i}(y))
\leftrightarrow\mathcal{D}_i(x)||\mathcal{D}_{i}(x)| = \frac
{1}{c_0}  \Biggl(c_0 \rho\prod_{j=1}^{i} [L(l_j)]^{2^{-j}} (l_0)^d
 \Biggr)^{2^i}  \Biggr)\\
&&\qquad \leq \mathbb{P} \Biggl(\Lambda_{i}(\bar{n}_{i}(y))
\leftrightarrow\mathcal{D}_i(x)||\mathcal{D}_{i}(x)| = \frac
{1}{c_0}  \Biggl(c_0 \rho\prod_{j=1}^{i} [L(l_j)]^{2^{-j}} (l_0)^d
 \Biggr)^{2^i}  \Biggr).
\end{eqnarray*}
By Lemma \ref{nontrivpcor} and $l_0> K_3(1/25)$, we see that
%
%e25 ###
%
\begin{eqnarray}\label{lemmacomp2}
&&\mathbb{P} \bigl(\mbox{$y$ is good up to level $0$}| \mbox{$y$ is
chosen uniformly}\nonumber\\[-8pt]\\[-8pt]
&&\qquad\hspace*{71.6pt}\mbox{at random from $\Lambda_{i}(\bar{n}_{i}(y))$} \bigr)
\geq\tfrac{24}{25} \rho.\nonumber
\end{eqnarray}
Furthermore,
\begin{eqnarray*}
&& \mathbb{P}^* \Biggl(\Lambda_{i}(\bar{n}_{i}(y))
\leftrightarrow\mathcal{D}_i(x)||\mathcal{D}_{i}(x)| = \frac
{1}{c_0}  \Biggl(c_0 \rho\prod_{j=1}^{i} [L(l_j)]^{2^{-j}} (l_0)^d
 \Biggr)^{2^i}  \Biggr) \\
&&\qquad =  1-\exp\Biggl[- \frac{1}{c_0} \Biggl(c_0 \rho\prod_{j=1}^{i}
[L(l_j)]^{2^{-j}} (l_0)^d  \Biggr)^{2^i}(l_0)^{d2^i}\lambda(l_{i+1})\Biggr]\\
&&\qquad =  1-\exp\Biggl[- \frac{1}{c_0}  \Biggl(c_0 \rho\prod_{j=1}^{i}
[L(l_j)]^{2^{-j}}  \Biggr)^{2^i}L (l_{i+1})\Biggr]\\
&&\qquad \geq \frac{L (l_{i+1})}{2c_0}  \Biggl(c_0 \rho\prod_{j=1}^{i}
[L(l_j)]^{2^{-j}}  \Biggr)^{2^i},
\end{eqnarray*}
where we have used the fact that $1-e^{-x} \geq x/2$ for $0< x \leq
1$, and the facts that $\rho, c_0<1$ and $l_0>K_2$.

Observe that $A_{i+1}(x)$ is dominated by a random variable which is
binomially distributed with parameters $n_i$ and $p_i$, where
\[
n_i = (l_{i+1}/l_i)^d - 1 = (l_0)^{d2^i}-1 \geq(l_0)^{d2^i}/2
\]
and
\[
p_i > \tfrac{24}{25} \rho\tfrac{2}{3} \bigl(1-\exp[-(l_0)^{d 2^i}
|\mathcal{D}_i(x)| \lambda(l_{i+1})]\bigr)
\]
by (\ref{lemmacomp}) and (\ref{lemmacomp2}).

If $x$ is good up to level $i$, then by (\ref{binafs2}), it holds that
\begin{eqnarray*}
p_i &>& \Biggl(\frac{16 \rho}{25} \Biggl(1-\exp\Biggl[-(l_0)^{d 2^i} \frac{1}{c_0}
\Biggl(c_0 \rho\prod_{j=1}^{i} [L(l_j)]^{2^{-j}} (l_0)^d
\Biggr)^{2^i} \lambda(l_{i+1})\Biggr]\Biggr)\Biggr)\\
&\geq& \frac{8 \rho L (l_{i+1})}{25 c_0}  \Biggl(c_0 \rho\prod
_{j=1}^{i} [L(l_j)]^{2^{-j}}  \Biggr)^{2^i}\\
&=& \frac{8 \rho}{25 c_0}  \Biggl(c_0 \rho[L
(l_{i+1})]^{2^{-(i+1)}}\prod_{j=1}^{i+1} [L(l_j)]^{2^{-j}}
\Biggr)^{2^i}\\
&\geq& \frac{8 \rho}{25 c_0}  \Biggl(c_0 \rho\Biggl[\prod
_{j=1}^{\infty} [L(l_j)]^{2^{-j}}\Biggr]^2  \Biggr)^{2^i}.
\end{eqnarray*}
The second line above, together with $c_0 = 2 \rho/25$, implies that
\[
m_{i+1} = L(l_{i+1})\Biggl(c_0 \rho(l_0)^d \prod_{j=1}^{i}
[L(l_j)]^{2^{-j}}\Biggr)^{2^i} = \frac{25 c_0}{4 \rho}\frac
{(l_0)^{d2^i}}{(l_0)^{d2^i}-1} \frac{n_ip_i}{2} < \frac{n_ip_i}{2}.
\]

By
%
%e26 ###
%
\begin{eqnarray}
\frac{n_ip_i}{2} &\geq& \frac{1}{4} (l_0)^{d 2^i} \frac{8 \rho}{25
c_0}  \Biggl(c_0 \rho\Biggl[\prod_{j=1}^{\infty} [L(l_j)]^{2^{-j}}\Biggr]^2  \Biggr)^{2^i}
\nonumber\\[-8pt]\\[-8pt]
&=& \frac{2 \rho}{25 c_0}  \Biggl(c_0 \rho\Biggl[\prod_{j=1}^{\infty}
[L(l_j)]^{2^{-j}}\Biggr]^2 (l_0)^d \Biggr)^{2^i},\nonumber
\end{eqnarray}
together with $l_0 > (100/\rho^2)^{1/(d-2/5)}$, $c_0 := 2 \rho/25$
and (\ref{binafs}), we obtain,
\begin{eqnarray*}
&&\mathbb{P}\bigl(A_{i+1}(x) \geq m_{i+1}| \mbox{$x$ is good up to level
$i$}\bigr) \\
&&\qquad \geq \mathbb{P}\biggl(A_{i+1}(x) \geq\frac{n_ip_i}{2}\bigg| \mbox
{$x$ is good up to level $i$}\biggr)\\
&&\qquad \geq 1 - \frac{4}{n_ip_i}\\
&&\qquad \geq 1 - \frac{25 c_0}{\rho}  \Biggl(c_0 \rho\Biggl[\prod
_{j=1}^{\infty} [L(l_j)]^{2^{-j}}\Biggr]^2 (l_0)^d \Biggr)^{-2^i}\\
&&\qquad \geq 1 -  \Biggl(\rho^2/25 \Biggl[\prod_{j=1}^{\infty}
[L(l_j)]^{2^{-j}}\Biggr]^2 (l_0)^d \Biggr)^{-2^i}\\
&&\qquad \geq 1 - \frac{1}{4^{2^i}}.
\end{eqnarray*}
This completes the proof of Lemma \ref{hierarchlem1}.
\end{pf}
\begin{pf*}{Proof of Lemma \ref{kolmouse}} The first step in the proof
is the observation that the event
\[
\mathcal{E} := \{\mbox{the number of ultimately good vertices in
$\mathbb{Z}^d$ is infinite}\}
\]
is independent of any finite set of edges. Indeed, for every finite set
of edges $E_0$, there is an $i \in\mathbb{N}$ such that all edges in
$E_0$ are shared by vertices in $\Lambda_i(0)$. However, whether a
vertex $x$ with $r(0,x)>i$ is ultimately good does not depend on edges
with at least one end-vertex in $\Lambda_i(0)$. Therefore, $\mathcal
{E}$ does not depend on $E_0$.

By a Kolmogorov-like zero--one law (see, e.g., \cite{Grim92},
page 289), we know that the probability that there will be
infinitely many ultimately good vertices is either 0 or 1. We will
prove that, with positive probability, every annulus of the form
$\Lambda_{i+1}(0) \setminus\Lambda_i(0)$ with $i \in\mathbb{N}$
contains at least one ultimately good vertex. This will prove the lemma.

Note that, by Lemma \ref{hierarchlem1},
\[
\mathbb{P}(\mbox{$x$ is ultimately good$| x$ is good up to level 0})
\geq1 - \sum_{i=0}^{\infty}\frac{1}{4^{2^i}} \geq2/3
\]
and, thus, by Lemma \ref{nontrivpcor},
\begin{eqnarray*}
&& \mathbb{P}(\mbox{$x$ is ultimately good}|x \in\mathcal
{C}_{\infty})\\
&&\qquad \geq \mathbb{P}(\mbox{$x$ is ultimately good$|x$ is good up to
level 0})\\
&&\qquad\quad{} \times\mathbb{P}(\mbox{$x$ is good up to level 0}|x \in
\mathcal{C}_{\infty})\\
&&\qquad \geq \tfrac{2}{3}\tfrac{24}{25}.
\end{eqnarray*}

The probability that the annulus $\Lambda_{i+1}(0) \setminus\Lambda
_i(0)$ contains no vertex that is good up to level $i$ is given by
\[
\bigl(\mathbb{P}\bigl(\mbox{$\Lambda_i(0) $ contains no vertex that is good up
to level $i$}\bigr)\bigr)^{(l_{i+1}/l_i) -1},
\]
where we have used the fact that the events that vertices in different
level $i$ blocks are good up to level $i$ are independent.
Note that by Lemma \ref{nontrivpcor} and $L_0>L_3(1/25)$, we have
\begin{eqnarray*}
&& \mathbb{P}\bigl(\mbox{$\Lambda_i(0)$ contains no vertex that is
good up to level $i$}\bigr)\\
&&\qquad \leq 1 - \mathbb{P}(\mbox{$x$ is ultimately good$|x$ is good up
to level $0$})\\
&&\qquad\quad{} \times\mathbb{P}\bigl(\mbox{$\Lambda_0(0)$ contains at least one
vertex that is good up to level $0$}\bigr)\\
&&\qquad \leq 1 - 48/75\\
&&\qquad = 9/25.
\end{eqnarray*}
Therefore, the probability that the annulus $\Lambda_{i+1}(0)
\setminus\Lambda_i(0)$ contains no vertex that is good up to level
$i$ is less than or equal to $(9/25)^{(d(l_0)^{2^i} -1)}$, which, in
turn, is less than $e^{-(16d/50)(l_0)^{2^i}}$, by $1-x \leq e^{-x}$.
Furthermore, by Lemma \ref{hierarchlem1}, it holds that
%
%e27 ###
%
\begin{eqnarray}
&&\mathbb{P}(\mbox{$x$ is ultimately good$|x$ is good up to level
$i$})\nonumber\\
&&\qquad\geq1 - \sum_{j=i}^{\infty} 4^{-2^j}
\geq1 - \sum_{j=i+1}^{\infty} 4^{-j} \\
&&\qquad\geq1-
(1/3)4^{-i}.\nonumber
\end{eqnarray}

For every $i\in\mathbb{N}$, the event that the annulus $\Lambda
_{i+1}(0) \setminus\Lambda_i(0) $ contains at least one ultimately
good vertex is increasing (for a definition of increasing events, see
\cite{Grim99}, page 32), so by the FKG inequality \cite
{Fort71,Grim99}, we obtain
\begin{eqnarray*}
&& \mathbb{P}\biggl( \bigcap_{i \in\mathbb
{N}}\bigl\{\bigl(\Lambda_{i+1}(0) \setminus\Lambda_i(0)\bigr)\mbox{ contains at
least one ultimately good vertex}\bigr\}\biggr)\\
&&\qquad \geq \prod_{i \in\mathbb{N}} \mathbb{P}\bigl(\mbox
{$\Lambda_{i+1}(0) \setminus\Lambda_i(0)$ contains at least one
ultimately good vertex}\bigr) \\
&&\qquad \geq \prod_{i \in\mathbb{N}}
\bigl(1-e^{-(16d/50)(l_0)^{2^i}}\bigr) \bigl(1- (1/3)4^{-i}\bigr)\\
&&\qquad \geq \prod_{i \in\mathbb{N}}
\bigl(1-e^{-(16d/50)l_02^i}\bigr) \bigl(1- (1/3)4^{-i}\bigr)\\
&&\qquad \geq 1- \sum_{i \in\mathbb{N}}
e^{-(16d/50)l_02^i} - \sum_{i \in\mathbb{N}} (1/3)4^{-i}\\
&&\qquad \geq 1- \sum_{i \in\mathbb{N}_+}
e^{-(16d/50)l_0i} - 4/9\\
&&\qquad \geq 1-\frac{e^{-(16d/50)l_0}}{1-e^{-(16d/50)l_0}}
- 4/9\\
&&\qquad \geq 1/18,
\end{eqnarray*}
where we have used the facts that $l_0 > 100/(16d)$ and
$e^{-x}(1-e^{-x})^{-1} < x^{-1}$ for $x>0$ in the
final inequality.
\end{pf*}

\section*{Acknowledgments}
The author wishes to thank Artem Sapozhnikov and Ronald Meester for
helpful discussions.

\printaddresses

\end{document}